\input amstex
\advance\vsize-0.5cm\voffset=-0.5cm\advance\hsize1cm\hoffset0cm
\magnification=\magstep1
\documentstyle{amsppt}
\NoBlackBoxes
\hfuzz=4.7pt
\hsize=6.5truein
\vsize=9truein
\topmatter
\bigskip
\title
{Splitting along a submanifold pair}
\endtitle
\author Rolando Jimenez, Yu.V. Muranov, and Du\v san Repov\v s
\endauthor
\thanks{Partially supported by  CONACyT and  DGAPA-UNAM,  by
the Russian Foundation for Fundamental
Research grant 05--01--00993,
and by the Ministry of Higher Education, Science and Technology
of the Republic of Slovenia research program P1--292--0101--04}
\endthanks
\endkeywords
\abstract{The paper
introduces a group $LSP$ of obstructions for  splitting  a homotopy
equivalence along a pair of submanifolds.
We develop exact sequences relating the $LSP$-groups with various surgery obstruction
groups for manifold triple and structure sets arising from
triples of manifolds. The natural map from the surgery obstruction group
of the ambient manifold to the $LSP$-group  provides
an invariant when elements of the Wall group are not realized by normal
maps of closed manifolds. Some $LSP$-groups are computed precisely.}
\endabstract
\endtopmatter
\document

\bigskip
\subhead 1. Introduction
\endsubhead
\bigskip

Consider a simple homotopy equivalence
$
f:M\to X
$
of closed $n$@-dimensional oriented topological manifolds.  Such a map is called
an {\it $s$@-triangulation of the manifold} $X$.
Two $s$-triangulations
$$
f_i:M_i\to X, \ i=1,2
$$
are said to be {\it equivalent}
if there exists an  orientation preserving
homeomorphism
$h:M_1\to M_2$
such that
the  diagram
$$
\matrix
M_1 &\overset{h}\to{\to} & M_2 \\
    &\searrow f_1  & \left\downarrow f_2\right. \\
    &                           & X
    \endmatrix
\tag 1.1
$$
is homotopy commutative.

The set of equivalence classes of  $s$@-triangulations of the manifold $X$
is denoted by $\Cal S(X)=\Cal S^{s}(X)$ (see \cite{18} and \cite{21}).
The computation  of the structure set $\Cal S^s(X)$ for a  manifold $X$
 is one the main
problems of geometric topology.

Let  $Y\subset X$ be a locally flat submanifold of codimension $q$
in $n$@-dimensional
topological manifold  $X$.
A simple homotopy equivalence $f: M\to X$ {\it splits along the submanifold} $Y$
(see \cite{18} and \cite{21}) if it is
homotopy equivalent to  a map $g$, transversal to $Y$  such that
    $N=g^{-1}(Y)$
satisfies the following properties:
$$
\matrix
i)& g|_N: N\to Y\  \text{is a simple homotopy equivalence},\\
ii) & g|_{(M \setminus N)}:M\setminus N\to X\setminus Y
\  \text{is a simple homotopy equivalence.}\\
\endmatrix
\tag 1.2
$$
A simple homotopy equivalence $g:M\to X$ with the properties (1.2)
is called  an $s$@-trian\-gulation of the pair $(X,Y)$. The set
of concordance classes of such $s$@-triangulations is denoted by
$\Cal S(X,Y, \xi)$ where $\xi$ is the topological normal
bundle of the submanifold $Y$ in $X$ (see \cite{18, \S 7.2}).

 Let $U$ be a tubular neighborhood of
the submanifold $Y$ in $X$, and let
 $\partial U$ denotes the boundary of $U$.
Denote by
$$
F=
\left(\matrix
\pi_1(\partial U)& \to &\pi_1(X\setminus Y)\\
\downarrow  & &     \downarrow\\
\pi_1(U)&\to &\pi_1(X)\\
\endmatrix\right)
\tag 1.3
$$
the push-out square of  fundamental groups with  orientations.

An obstruction to splitting  the map $f$ along the submanifold
$Y$
lies in the splitting obstruction group
$LS_{n-q}(F)$  which
depends only on $n-q \bmod 4 $ and on the push-out square $F$.

In  fact, the  obstruction to splitting  defines correctly the
map \cite{18} that fits  in
the following exact sequence
$$
\cdots \to \Cal S(X,Y, \xi)\to \Cal S(X)\to LS_{n-q}(F).
\tag 1.4
$$

The splitting obstruction groups are closely related
to other obstruction groups which arise naturally
for the manifold pair $Y\subset X$ (see \cite{1}, \cite{2}, \cite{13}, \cite{18}, and \cite{21}).
The main relation is
given by the following braid of exact sequences
(see \cite{18} and \cite{21})
$$
\smallmatrix
\rightarrow & L_n(\pi_1(X\setminus Y) & \longrightarrow &
L_{n}(\pi_1(X)) &
\longrightarrow & {LS}_{n-q-1}(F) & \rightarrow\cr
\ &  \nearrow \ \ \ \ \ \ \ \ \searrow & \ &  \nearrow \ \ \ \ \ \ \
\
\searrow
& \  & \nearrow \ \ \ \ \ \ \ \  \searrow & \ \cr
\ & \ & LP_{n-q}(F)& \ & L_{n}(\pi_1(X\setminus Y)\to\pi_1(X)) & \ & \ \cr
\ &  \searrow \ \ \ \ \ \ \ \ \nearrow & \ &  \searrow \ \ \ \ \ \ \
\
\nearrow
&   & \searrow \ \ \ \ \ \ \ \  \nearrow & \ \cr
\rightarrow & {LS}_{n-q}(F) & \longrightarrow &
L_{n-q}(\pi_1(Y)) &
\rightarrow & L_{n-1}(\pi_1(X\setminus Y))& \rightarrow
\endmatrix
\tag 1.5
$$
where $L_*=L^s_*$ denote  the surgery obstruction groups and
$LP_*(F)=LP_*^s(F)$ denote the surgery obstruction groups of
the manifold pair $(X, Y)$.  The groups $LP_*(F)$ also
depend  only on $n-q\bmod 4$ and on the square $F$.

The main methods for computing  the set
$\Cal S(X)$ (for $n\geq 4$) are based on  the surgery exact sequence
(see \cite{17}, \cite{18}, and \cite{21})
$$
\cdots \to L_{n+1}(\pi_1(X))\to \Cal S(X)\to [X, G/TOP]\overset{\sigma}\to
{\to} L_n(\pi_1(X))
\tag 1.6
$$
where the set $[X, G/TOP]$ is isomorphic to the set of
concordance classes of topological normal maps to the manifold $X$.

The set $\Cal S(X,Y, \xi)$ fits into the surgery exact sequence \cite{18, page 584}
for the manifold pair $(X,Y)$
$$
\cdots \to LP_{n-q+1}(F)\to \Cal S(X,Y, \xi)\to [X, G/TOP]\to
 LP_{n-q}(F).
\tag 1.7
$$
The exact sequence (1.7) is the natural  generalization
of  the exact sequence (1.6) to the case
of a manifold pair.

The computation of the map $\sigma$ in (1.6)
 is the basic step in investigating
the surgery exact sequence.  For manifolds with finite
fundamental groups  deep results in this direction
were obtained in
\cite{5}, \cite{6}, \cite{9}, \cite{10},  and \cite{11}.
The results of these papers are based on relations between
the surgery exact sequence and the splitting problem
for a one-sided submanifold.

Let
$$
Z^{n-q-q^{\prime}}\subset Y^{n-q}\subset X^n
\tag 1.8
$$
be a triple of
closed topological manifolds.
We shall consider only
locally flat topological submanifolds equipped
with the structure of a normal topological bundle (see
\cite{18, pages 562--563}). Such a triple of manifolds
defines a stratified manifold $\Cal X$ in the sense of Browder and Quinn
(see \cite{4}, \cite{14}, \cite{15}, \cite{16},  and \cite{22}).

A simple homotopy equivalence $f: M\to X$ is an {\it $s$@-triangulation of the
triple} if  every pair of manifolds from this triple
satisfies  properties that are similar to (1.2) for the pair $(X,Y)$
(see \cite{4}, \cite{16}, and \cite{21}).
The set of concordance classes of such $s$@-triangulations
is denoted by $\Cal S(\Cal X)=\Cal S(X, Y,Z)$.

Surgery theory is applicable to stratified spaces, and we have
the following
exact sequence (see \cite{4} and  \cite{22})
$$
\cdots \to L_{n+1}^{BQ}(\Cal X)\to \Cal S(\Cal X)\to [X, G/TOP]\to L_n^{BQ}(\Cal X)
\tag 1.9
$$
where $L_*^{BQ}(\Cal X)$ are the Browder-Quinn surgery obstruction groups
of the stratified space $\Cal X$. For these groups we
have isomorphisms
$$
L_n^{BQ}(\Cal X)=LT_{n-q-q^{\prime}}(X,Y,Z)
$$
with surgery obstruction groups $LT_*$ of the manifold triple $(X,Y,Z)$
(see \cite{14} and \cite{16}).

In the present paper we develop surgery theory for  manifold
triples in order  to investigate splitting a homotopy equivalence
along a submanifold pair.
By definition,  a simple homotopy equivalence
 $f:M\to X$ {\it splits along the submanifold pair} $(Z\subset Y)$
 if it is concordant  to an $s$@-triangulation $g$ of
 the triple $Z\subset Y\subset X$.
      We introduce  groups $LSP_*$  of obstructions to splitting a
simple homotopy equivalence $f: M\to X$ along a pair of embedded
submanifolds $(Z\subset Y)\subset X$ and describe their  relations to
classical obstruction groups in surgery theory. The group $LSP_*$
is a natural straightforward generalization of the group $LS_*$ if
we consider a pair of submanifolds $(Z\subset Y)$ instead of a
submanifold $Y$.
The $LSP$@-groups give in a natural way  an
invariant for determining when elements of Wall groups are not
realized by normal maps of closed manifolds.
This invariant is equivalent to the pair of Hambleton's  invariants
($A$ and $B$) in paper \cite{6}.

The rest of the paper is organized as follows. In section 2, we recall
notation, constructions and results from the literature, which will be
needed in the current paper.
In section 3, we construct the spectrum $\Bbb L SP(X, Y Z)$ and
relate via exact sequences its homotopy groups $LSP_{\ast} (X, Y, Z)$
to classical obstruction groups and structure sets arising from
triples of manifolds. In section 4, we apply the above to obtain
results when elements of Wall groups are not realized by normal maps
of closed manifolds and compute some $LSP_{\ast}$-groups.

\bigskip
\subhead 2. Preliminaries
\endsubhead
\bigskip
The current paper will make significant use of constructions, ideas,
and results
in the papers \cite{1},\cite{2}, \cite{8}, \cite{13}, \cite{16}, \cite{17}, \cite{18}, and \cite{21}.
A thread running through all of these articles is the use, due to
Ranicki \cite{17}, \cite{18}, of spectra for developing the algebraic theory of surgery.
In this section we recall
some necessary definitions and results from these papers.

Consider  a triple of topological manifolds (1.8). Let
 $\xi$ denote the normal bundle of $Y$ in $X$ and $F$
 the square of fundamental groups
 in the splitting problem  for the pair $Y\subset X$.
Similarly
 we introduce the following bundles and squares:

 the bundle  $\eta$ and the square  $\Psi$
for the pair $Z\subset Y$,

the bundle   $\nu$ and the square  $\Phi$
 for the pair $Z\subset X$.

Let $U_{\xi}$ be the space of the normal bundle $\xi$.
We shall assume that the space $U_{\nu}$ of the
normal bundle  $\nu$ is identified with the space
$V_{\xi}$ of the restriction of the bundle $\xi$ to the space
$U_{\eta}$ of the normal bundle $\eta$ so that
$\partial U_{\nu}=\partial U_{\xi}|_{U_{\eta}}\cup
U_{\xi}|_{\partial U_{\eta}}$ (see \cite{4}, \cite{15}, \cite{16},
 and \cite{22}).

The conditions on the spaces of normal bundles for the manifold triple
(1.8) yield   a  pair of manifolds with boundaries
$$
(Y\setminus Z, \partial(Y\setminus Z))\subset (X\setminus Z,
\partial(X\setminus Z))
\tag 2.1
$$
where
$$
\partial(Y\setminus Z)\subset \partial(X\setminus Z)
\tag 2.2
$$
is
a closed manifold pair.
Denote by $F_Z$ the square of fundamental groups in the
splitting problem relative to
boundary for the  pair (2.1), and by $F_U$  the square in
the splitting problem for the pair (2.2).

For an arbitrary group $\pi$ with orientation the surgery $\Omega$-spectrum
$\Bbb L(\pi)=\Bbb L(\Bbb Z\pi)$ is defined (see \cite{8}, \cite{17}, and \cite{21}).
Here $\Bbb Z \pi$ denotes the integral group ring
equipped with the involution
$$
\Sigma a_g g \mapsto \Sigma a_g w(g)g^{-1}, \ a_g\in \Bbb Z, g\in \pi
$$
where $w:\pi \to\{\pm 1\}$ is the orientation homomorphism.
Recall that for this  $\Omega$@-spectrum we have
 $$
 \pi_n(\Bbb L(\pi))= L_n(\pi).
 $$
Let $\bold L_{\bullet}$ denote the 1-connected cover of the spectrum
$\Bbb L(1)$ with
${\bold L_{\bullet}}_0=G/TOP$.
For a topological space $X$ we have the following
cofibration (see \cite{17} and \cite{18})
$$
X_+\land \bold L_{\bullet}\to \Bbb L(\pi_1(X))\to \Bbb S(X) .
\tag 2.3
$$
The homotopy long exact sequence of the cofibration (2.3)
gives the algebraic surgery exact sequence of Ranicki \cite{17}
$$
\cdots \to L_{n+1}(\pi_1(X))\to \Cal S_{n+1}(X)\to H_n(X, \bold L_{\bullet})
\to  L_n(\pi_1(X))\to \cdots
\tag 2.4
$$
with
$$
\pi_{n+1}(\Bbb S(X))=\Cal S_{n+1}(X)\cong\Cal S^{TOP}(X).
$$
The left  part of the exact sequence (2.4) is isomorphic to the exact sequence
(1.6).

A similar result is valid  for the exact sequences
(1.4), (1.7), and (1.9).
In particular, we have cofibrations of spectra
$$
\Bbb S(X,Y, \xi)\to \Bbb S(X) \to \Sigma^{q+1} \Bbb LS(F),
\tag 2.5
$$
$$
X_+\land \bold L_{\bullet}\to \Sigma^q \Bbb LP(F)\to \Bbb S(X, Y, \xi),
\tag 2.6
$$
and
$$
X_+\land \bold L_{\bullet}\to \Sigma^{q+q^{\prime}}\Bbb
LT(X,Y,Z)\to \Bbb S(X,Y,Z),
\tag 2.7
$$
where $\Sigma$ denotes the suspension functor on the category of
$\Omega$@-spectra.
These cofibrations  generate exact sequences that contain parts
which are isomorphic to the  exact
 sequences (1.4), (1.7), and (1.9), respectively.

Recall that for an arbitrary pair
$(X,Y)$ of topological
spaces equipped with orientation, a spectrum
$\Bbb S(X,Y)$  for the relative structure sets
$\Cal S_*(X,Y)$ is defined (see \cite{17} and \cite{18}).

A homomorphism of oriented groups  $f : \pi\to \pi^{\prime}$
induces a cofibration of  $\Omega$@-spectra
$$
\CD
\Bbb L(\pi) @>>> \Bbb L(\pi') @>>> \Bbb L(f)
\endCD
\tag 2.8
$$
where  $\Bbb L(f)$ is the spectrum for relative $L$@-groups of
the map $f$.

For the  manifold pair    $(X,Y)$
we have a homotopy commutative diagram of spectra
(see \cite{1}, \cite{2}, \cite{8}, and \cite{18})
$$
\matrix
\Bbb L(\pi_1(Y)) & \overset{}\to{\rightarrow}
&\Sigma^{-q}\Bbb
L(\pi_1(\partial U)\to \pi_1(U))&
\overset{}\to{\rightarrow} & \Sigma^{-q}\Bbb
L(\pi_1(X\setminus Y)\to \pi_1(X))  \\
  & \ \searrow & \downarrow & & \downarrow \\
 & & \Sigma^{1-q}\Bbb L(\pi_1(\partial U))&
 \overset{}\to{\rightarrow} &\Sigma^{1-q}\Bbb
L(\pi_1(X\setminus Y)),
\endmatrix
\tag 2.9
$$
where the left maps are transfer maps on the spectra level,
and  the right horizontal maps are induced by inclusions.

The diagram (2.9)  provides a homotopy commutative diagram of spectra
$$
\matrix
\Bbb L(\pi_1(Y)) & \overset{}\to{\rightarrow} & \Sigma^{-q}\Bbb L(\pi_1(X\setminus Y)\to \pi_1(X))&\to& \Sigma\Bbb LS(F)  \\
 \downarrow= &  & \downarrow & & \downarrow \\
\Bbb L(\pi_1(Y)) & \overset{}\to{\rightarrow} &\Sigma^{1-q}\Bbb L(\pi_1(X\setminus Y))& \to& \Sigma\Bbb LP(F)\\
\endmatrix
\tag 2.10
$$
in which the horizontal rows  are cofibrations.
The  homotopy long exact sequences of the maps from diagram  (2.10)
generate the diagram (1.5).

The triple (1.8) defines also on the spectra level
the  maps  (see \cite{16} and \cite{22})
$$
\Bbb L(\pi_1(Z))\to \Sigma^{-q^{\prime}+1}\Bbb LP(F_U)\to
\Sigma^{-q^{\prime}+1}\Bbb LP(F_Z)
\tag 2.11
$$
where the first map is the transfer map, and the second
map is induced by the inclusion in (2.1).

By \cite{16} and \cite{22} we have
a cofibration
$$
\Bbb L(\pi_1(Z))\to
\Sigma^{-q^{\prime}+1}\Bbb LP(F_Z)\to \Sigma\Bbb LT(X,Y,Z)
\tag 2.12
$$
where the first map is the composition of the maps in (2.11).

Consider the composition of the maps
 $$
\Bbb LP(F)\to \Bbb L(\pi_1(Y))\to \Bbb S(Y)\to \Sigma^{q^{\prime}+1}
\Bbb LS(\Psi).
\tag 2.13
$$
The first map in (2.13) follows from (2.10), the second is the map
from (1.3) for the manifold $Y$, and the third map is the map
from (2.5) for the pair $(Y,Z)$.
By \cite{14} and \cite{16} we have the
 cofibration
$$
\Bbb LP(F)\to \Sigma^{q^{\prime}+1}\Bbb LS(\Psi) \to
\Sigma^{q^{\prime}+1}\Bbb LT(X,Y,Z).
\tag 2.14
$$
From the cofibration (2.14), we obtain the  homotopy
pull-back square of spectra
$$
\matrix
\Bbb LT(X,Y,Z)&\to & \Sigma^{-q^{\prime}}\Bbb LP(F) \\
\downarrow &&\downarrow \\
\Bbb LP(\Psi)&\to & \Sigma^{-q^{\prime}}\Bbb L(\pi_1(Y)) \\
\endmatrix
\tag 2.15
$$
where  the cofibres of the vertical maps are naturally homotopy equivalent
to the spectrum $\Sigma^{-q-q^{\prime}+1}\Bbb L(\pi_1(X\setminus Y))$.

Consider the commutative diagram of inclusions
$$
\matrix
(Y\setminus Z)&\subset &(X\setminus Z)\\
\cap& & \cap\\
Y&\subset & X.\\
\endmatrix
\tag 2.16
$$
The horizontal inclusions of submanifolds of codimension $q$, provide
as in (2.10),   the transfer maps   fitting
into the homotopy commutative diagram
$$
\matrix
\Bbb L(\pi_1(Y\setminus Z)) & \overset{}\to{\rightarrow} & \Sigma^{-q}\Bbb L(\pi_1(X\setminus Y)\to \pi_1(X\setminus Z))&\to& \Sigma\Bbb LS(F_Z)  \\
 \downarrow &  & \downarrow & & \downarrow \\
\Bbb L(\pi_1(Y)) & \overset{}\to{\rightarrow} &\Sigma^{-q}\Bbb L(\pi_1(X\setminus Y)\to X)& \to& \Sigma\Bbb LS(F)\\
 \downarrow &  & \downarrow & & \downarrow \\
\Bbb L(\pi_1(Y\setminus Z)\to (\pi_1(Y)) & \overset{tr^{rel}}\to{\rightarrow} &\Sigma^{-q}\Bbb L(\pi_1(X\setminus Z)\to \pi_1(X))& \to& \Sigma^{1+q^{\prime}}\Bbb LNS\\
\endmatrix
\tag 2.17
$$
in which  the upper vertical maps are induced by the vertical maps from  (2.16).
The spectrum $\Bbb LNS=\Bbb LNS(X,Y,Z)$ is the
spectrum for the relative $L$@-groups
of the  map $tr^{rel}$
(see \cite{7} and \cite{15})
with the homotopy groups
$$
LNS_n =LNS_n(X,Y,Z)=\pi_n(\Bbb LNS).
$$
Note that the diagram (2.17) generates the following commutative diagram
\cite{15}
$$
\smallmatrix
&\vdots     &  &\vdots     & & \vdots   & \\
&\downarrow &  &
\downarrow &
 &\downarrow & \\
\dots\rightarrow &LS_{n-q}(F_Z)&\rightarrow &
L_{n-q}(\pi_1(Y\setminus Z))&
\overset{}\to{\rightarrow} &
L_{n}(\pi_1(X\setminus Y)\to \pi_1(W)) &\rightarrow \cdots \\
&\downarrow &  &
\downarrow &
 &\downarrow & \\
\dots\rightarrow & LS_{n-q}(F)&\rightarrow &
L_{n-q}(\pi_1(Y))&
\overset{}\to{\rightarrow} &L_{n}(\pi_1(X\setminus Y)\to \pi_1(X)) &\rightarrow \cdots\\
&\downarrow &  &\downarrow &&\downarrow & \\
\dots\rightarrow & LNS_{k}&\rightarrow &
L_{n-q}(\pi_1(Y\setminus Z)\to \pi_1(Y))&
\overset{}\to{\rightarrow}
&L_{n}(\pi_1(W)\to \pi_1(X)) &\rightarrow \cdots\\
 &\downarrow &  &\downarrow & &\downarrow & \\
 &\vdots     &  &\vdots     & & \vdots   & \\
\endmatrix
\tag 2.18
 $$
where $k=n-q-q^{\prime}$ and $W=X\setminus Z$.

\bigskip
\subhead 3. Splitting a homotopy equivalence along a submanifold pair
\endsubhead
\bigskip

For the  triple of  manifolds (1.8)
we introduce  below the spectrum $\Bbb  LSP(X,Y,Z)$
with homotopy
groups
$$
LSP_*= LSP_*(X,Y,Z)=\pi_n(\Bbb LSP(X,Y,Z)).
\tag 3.1
$$
The groups $LSP_*(X,Y,Z)$ are a natural straightforward
generalization of the splitting obstruction groups $LS_*(F)$
to the case when the
manifold $X$ contains a pair of embedded submanifolds $(Z\subset
Y)\subset X$ instead of first a single submanifold $Y$.
We describe via exact sequences the relation of  the groups $LSP_*(X,Y,Z)$
to classical obstruction groups and structure sets which arise naturally
for a triple of manifolds.

The bottom map in the diagram (2.15) and the commutative
diagram (2.10)
provide the homotopy  commutative diagram of spectra
$$
\matrix
\Bbb LP(\Psi) & \overset{}\to{\rightarrow} & \Sigma^{-q-q^{\prime}}\Bbb L(\pi_1(X\setminus Y)\to \pi_1(X)) \\
 \downarrow= &  & \downarrow   \\
\Bbb LP(\Psi) & \overset{}\to{\rightarrow} &\Sigma^{1-q-q^{\prime}}\Bbb L(\pi_1(X\setminus Y))\\
\endmatrix
\tag 3.2
$$
in which the fiber of the bottom map is the spectrum $\Bbb LT(X,Y,Z)$.
This follows from the pull-back property of the square (2.15).
Denote by $\Bbb LSP(X,Y,Z)$  the fiber of the upper horizontal map in (3.2).
We obtain the homotopy commutative diagram of spectra
$$
\matrix
\Bbb LP(\Psi) & \overset{}\to{\rightarrow} & \Sigma^{-q-q^{\prime}}\Bbb L(\pi_1(X\setminus Y)\to \pi_1(X))&\to& \Sigma\Bbb LSP(X,Y,Z)\\
 \downarrow= &  & \downarrow  && \downarrow  \\
\Bbb LP(\Psi) & \overset{}\to{\rightarrow} &\Sigma^{1-q-q^{\prime}}\Bbb L(\pi_1(X\setminus Y))&\to&\Sigma  \Bbb LT(X,Y,Z)\\
\endmatrix
\tag 3.3
$$
in which the right vertical map is induced by the  two others vertical maps
(see \cite{20}). Note  that   the right square in (3.3) is a pull-back.

\proclaim{Proposition 3.1} The groups $LSP_*(X,Y,Z)$ that are defined by
(3.1) fit into the following braid of exact sequences
$$
\matrix
\rightarrow & L_{n}(C) & \longrightarrow &
L_n( \pi_1(X)) &
\rightarrow & LSP_{k-1}& \rightarrow \cr
\ &  \nearrow \ \ \ \ \ \ \ \ \searrow & \ &  \nearrow \ \ \ \ \ \ \
\
\searrow
& \  & \nearrow \ \ \ \ \ \ \ \  \searrow & \ \cr
\ & \ & LT_k(X,Y,Z)& \ & L_{n}(C\to D) & \ & \ \cr
\ &  \searrow \ \ \ \ \ \ \ \ \nearrow & \ &  \searrow \ \ \ \ \ \ \
\
\nearrow
& \  & \searrow \ \ \ \ \ \ \ \  \nearrow & \ \cr
\rightarrow & LSP_{k} & \longrightarrow &
LP_{k}(\Psi) &
\longrightarrow & L_{n-1}(C) & \rightarrow,
 \endmatrix
\tag 3.4
$$
where $C=\pi_1(X\setminus Y)$, $D= \pi_1(X)$, and  $k=n-q-q^{\prime}$.
The diagram (3.4) is realized on the spectra level.
\endproclaim
\demo{Proof} The right square in the diagram (3.3) is  a pull-back.
The
homotopy long exact sequences of this square provide the
commutative braid of exact sequences (3.4).
\qed
\enddemo
\bigskip

\proclaim{Theorem 3.2} There exists a commutative braid  of exact
sequences
$$
\matrix
\rightarrow & \Cal S_{n+1}(X,Y,Z) & \longrightarrow &
H_n(X,\bold  L_{\bullet}) &
\rightarrow & L_{n}(\pi_1(X))& \rightarrow \cr
\ &  \nearrow \ \ \ \ \ \ \ \ \searrow & \ &  \nearrow \ \ \ \ \ \ \
\
\searrow
& \  & \nearrow \ \ \ \ \ \ \ \  \searrow & \ \cr
\ & \ & \Cal S_{n+1}(X)& \ & LT_{n-q-q^{\prime}} & \ & \ \cr
\ &  \searrow \ \ \ \ \ \ \ \ \nearrow & \ &  \searrow^{\alpha} \ \ \ \ \ \ \
\
\nearrow
& \  & \searrow \ \ \ \ \ \ \ \  \nearrow & \ \cr
\rightarrow & L_{n+1}(\pi_1(X)) & \longrightarrow &
LSP_{n-q-q^{\prime}} &
\longrightarrow & \Cal S_{n}(X,Y,Z) & \rightarrow
 \endmatrix
\tag 3.5
$$
which is realized on the spectra level.
\endproclaim

\demo{Proof} Consider the homotopy commutative square of spectra
$$
\matrix
X_+\land \bold  L_{\bullet} &\to & \Bbb L(\pi_1(X)) \\
\downarrow & & \downarrow=\\
\Sigma^{q+q^{\prime}}\Bbb LT& \to & \Bbb L(\pi_1(X))\\
\endmatrix
\tag 3.6
$$
in which the upper horizontal map lies in (2.3), the left vertical map
lies in (2.7), and the bottom horizontal map is the map from the
diagram (3.4) on spectra level (see \cite{14} and \cite{15}).
The diagram (3.6) induces a map of the fibres of its horizontal maps.
We obtain the homotopy commutative diagram of spectra
$$
\matrix
\Sigma^{-1}\Bbb S(X)&\to&X_+\land \bold  L_{\bullet} &\to & \Bbb L(\pi_1(X)) \\
\downarrow& &\downarrow & & \downarrow=\\
\Sigma^{q+q^{\prime}}\Bbb LSP&\to&
\Sigma^{q+q^{\prime}}\Bbb LT& \to & \Bbb L(\pi_1(X))\\
\endmatrix
$$
in which the left square is a push-out.
The homotopy long exact sequences of this square give
the diagram (3.5).
\qed
\enddemo
\smallskip

The commutative diagram (3.5) is a natural generalization of
the  diagram in \cite{18, Proposition 7.2.6, iv} to the case
of a pair of submanifolds $Z\subset Y$ in the manifold $X$.
The map
$$
\Sigma^{-1}\Bbb S(X)\to\Sigma^{q+q^{\prime}}\Bbb LSP
$$
 induces a map
$$
\alpha:\Cal
S_{n+1}(X)\to LSP_{n-q-q^{\prime}}(X,Y,Z)
$$
that on the algebraic level
corresponds to taking  the obstruction to splitting
along the submanifold pair $Z\subset Y$.

Now we describe  the relation of $LSP_*$ to classical surgery obstruction
groups for the triple $(X,Y,Z)$  of manifolds    (1.8).
\smallskip

\proclaim{Theorem 3.3} There exist braids of exact sequences
$$
\matrix
\rightarrow & LS_{n-q}(F_Z) & \longrightarrow &
LT_k(X,Y,Z) &
\rightarrow & L_{n}(\pi_1(X))& \rightarrow \cr
\ &  \nearrow \ \ \ \ \ \ \ \ \searrow & \ &  \nearrow \ \ \ \ \ \ \
\
\searrow
& \  & \nearrow \ \ \ \ \ \ \ \  \searrow & \ \cr
\ & \ & LSP_k& \ & LP_{k}(\Phi) & \ & \ \cr
\ &  \searrow \ \ \ \ \ \ \ \ \nearrow & \ &  \searrow \ \ \ \ \ \ \
\
\nearrow
& \  & \searrow \ \ \ \ \ \ \ \  \nearrow & \ \cr
\rightarrow & L_{n+1}(\pi_1(X)) & \longrightarrow &
LS_{k}(\Phi) &
\longrightarrow & LS_{n-q-1}(F_Z) & \rightarrow,
 \endmatrix
\tag 3.7
$$
$$
\matrix
\rightarrow & LS_{n-q}(F_Z) & \longrightarrow &
LS_{n-q}(F) &
\rightarrow & LS_{k-1}(\Psi)& \rightarrow \cr
\ &  \nearrow \ \ \ \ \ \ \ \ \searrow & \ &  \nearrow \ \ \ \ \ \ \
\
\searrow
& \  & \nearrow \ \ \ \ \ \ \ \  \searrow & \ \cr
\ & \ & LSP_k& \ & LNS_{k} & \ & \ \cr
\ &  \searrow \ \ \ \ \ \ \ \ \nearrow & \ &  \searrow \ \ \ \ \ \ \
\
\nearrow
& \  & \searrow \ \ \ \ \ \ \ \  \nearrow & \ \cr
\rightarrow & LS_{k}(\Psi) & \longrightarrow &
LS_{k}(\Phi) &
\longrightarrow & LS_{n-q-1}(F_Z) & \rightarrow,
 \endmatrix
\tag 3.8
$$
and
$$
\matrix
\rightarrow & LS_{n-q+1}(F) & \longrightarrow &
LS_k(\Psi) &
\rightarrow & LT_{k}& \rightarrow \cr
\ &  \nearrow \ \ \ \ \ \ \ \ \searrow & \ &  \nearrow \ \ \ \ \ \ \
\
\searrow
& \  & \nearrow \ \ \ \ \ \ \ \  \searrow & \ \cr
\ & \ & LP_{n-q+1}(F)& \ & LSP_{k} & \ & \ \cr
\ &  \searrow \ \ \ \ \ \ \ \ \nearrow & \ &  \searrow \ \ \ \ \ \ \
\
\nearrow
& \  & \searrow \ \ \ \ \ \ \ \  \nearrow & \ \cr
\rightarrow & LT_{k+1} & \longrightarrow &
L_{n+1}(\pi_1(X)) &
\longrightarrow & LS_{n-q}(F) & \rightarrow,
 \endmatrix
\tag 3.9
$$
where  $k=n-q-q^{\prime}$.
The braids (3.7), (3.8), and (3.9) are realized on the spectra level.
\endproclaim

\demo{Proof} The natural forgetful maps (see \cite{15} and \cite{16})
$$
\Bbb LT(X,Y,Z)\to \Bbb LP(\Phi)\to \Sigma^{-q-q^{\prime}}\Bbb L(\pi_1(X))
$$
provide the  homotopy commutative square
$$
\matrix
\Bbb LT(X,Y,Z) &\to &\Sigma^{-q-q^{\prime}}\Bbb L(\pi_1(X)) \\
\downarrow & &\downarrow= \\
 \Bbb LP(\Phi)&\to&
\Sigma^{-q-q^{\prime}}\Bbb L(\pi_1(X)).\\
\endmatrix
\tag 3.10
$$
The square induces a map of the fibres of its horizontal maps (see \cite{20}).
Thus we obtain a homotopy commutative diagram
$$
\matrix
\Bbb LSP(X,Y,Z) &\to&\Bbb LT(X,Y,Z) &\to &
\Sigma^{-q-q^{\prime}}\Bbb L(\pi_1(X)) \\
\downarrow & &\downarrow & &\downarrow= \\
\Bbb LS(\Phi) & \to & \Bbb LP(\Phi)&\to&
\Sigma^{-q-q^{\prime}}\Bbb L(\pi_1(X))\\
\endmatrix
$$
in which the left square is a push-out (and hence a pull-back).
Now, similarly to
Proposition 3.1, we obtain the  diagram (3.7).

The diagram (2.10) for the pair $(Y,Z)$ provides a homotopy
commutative pull-back square of spectra
$$
\matrix
\Bbb LP(\Psi) &\to & \Sigma^{-q^{\prime}}\Bbb L(\pi_1(Y))\\
\downarrow & &\downarrow \\
\Bbb L(\pi_1(Z)) &\to& \Sigma^{-q^{\prime}}\Bbb L(\pi_1(Y\setminus Z) \to \pi_1(Y)). \\
\endmatrix
\tag 3.11
$$
There is a homotopy commutative pull-back square of spectra
$$
\matrix
\Sigma^{-q^{\prime}-q}\Bbb L(\pi_1(X\setminus Y)\to \pi_1(X))
&\overset{=}\to{\to} &
\Sigma^{-q^{\prime}-q}\Bbb L(\pi_1(X\setminus Y)\to \pi_1(X))\\
\downarrow & &\downarrow \\
\Sigma^{-q^{\prime}-q}\Bbb L(\pi_1(X\setminus Z)\to \pi_1(X)) &
\overset{=}\to{\to}&
\Sigma^{-q^{\prime}-q}\Bbb L(\pi_1(X\setminus Z)\to \pi_1(X)) \\
\endmatrix
\tag 3.12
$$
in which the vertical maps are induced by the natural inclusion.
The transfer maps and diagrams (2.17) and (3.3) give the map
of diagram (3.11) to diagram (3.12). The cofibers of this
map of diagrams
provide a homotopy commutative pull-back square of spectra (see \cite{19})
$$
\matrix
\Sigma\Bbb LSP &\to & \Sigma^{-q^{\prime}+1}\Bbb LS(F)\\
\downarrow & & \downarrow \\
\Sigma\Bbb LS(\Phi) &\to &\Sigma \Bbb LNS. \\
\endmatrix
\tag 3.13
$$
This follows from (2.9), (2.18), and (3.3).
The diagram (3.8) follows from the square (3.13), similarly to
the previous case.

The natural forgetful maps in the diagram (2.15)
 $$
\Bbb  LT(X,Y,Z)\to
\Sigma^{-q^{\prime}}\Bbb LP(F) \to\Sigma^{-q-q^{\prime}}\Bbb L(\pi_1(X))
$$
provide
the homotopy commutative diagram of spectra
 $$
\matrix
\Bbb LT(X,Y,Z) &\to&\Sigma^{-q^{\prime}}\Bbb LP(F) &\to &
\Sigma\Bbb LS(\Psi) \\
\downarrow= & &\downarrow & &\downarrow \\
\Bbb LT(X,Y,Z) & \to & \Sigma^{-q-q^{\prime}}\Bbb L(\pi_1(X))&\to&
\Sigma\Bbb LSP(X,Y,Z),\\
\endmatrix
\tag 3.14
$$
in which the rows are cofibrations, and the right vertical map is
defined by \cite{19}. Hence  the right square in (3.14) is
a pull-back. From this, the diagram
(3.9) follows. \qed
\enddemo
\smallskip

\proclaim{Corollary 3.4} There exist exact sequences
$$
\cdots\to LSP_k \to LS_{n-q}(F)\to LS_{k-1}(\Psi)\to \cdots,
$$
$$
\cdots\to LSP_k \to LS_{k}(\Phi)\to LS_{n-q-1}(F_Z)\to \cdots,
$$
and
$$
\cdots\to LSP_k \to LP_{k}(\Psi)\to L_{n-1}(\pi_1(X\setminus Y)\to\pi_1(X))
\to \cdots,
$$
in which the  left maps are natural forgetful maps.
\endproclaim
\smallskip

Now we  describe some relations between the $LSP_*$-groups
 and  various structure sets which arise for the
triple of manifolds $(X,Y,Z)$.

\proclaim{Theorem 3.5} There exist  braids of exact sequences
$$
\matrix
\rightarrow & \Cal S_{n}(X) & \longrightarrow &
LSP_{k-1} &
\rightarrow &\Cal S_{l-1}(Y,Z,\eta)& \rightarrow \cr
\ &  \nearrow \ \ \ \ \ \ \ \ \searrow & \ &  \nearrow \ \ \ \ \ \ \
\
\searrow
& \  & \nearrow \ \ \ \ \ \ \ \  \searrow & \ \cr
\ & \ & \Cal S_{n}(X, X\setminus Y)& \ & \Cal S_{n-1}(X,Y,Z) & \ & \ \cr
\ &  \searrow \ \ \ \ \ \ \ \ \nearrow & \ &  \searrow \ \ \ \ \ \ \
\
\nearrow
& \  & \searrow \ \ \ \ \ \ \ \  \nearrow & \ \cr
\rightarrow & \Cal S_{l}(Y,Z,\eta) & \longrightarrow &
\Cal S_{n-1}(X\setminus Y) &
\longrightarrow & \Cal S_{n-1}(X) & \rightarrow,
 \endmatrix
\tag 3.15
$$

$$
\matrix
\rightarrow & H_{l}(Y, \bold L_{\bullet}) & \longrightarrow &
L_n(\pi_1(X\setminus Y)\to \pi_1(X)) &
\rightarrow & LSP_{k-1}& \rightarrow \cr
\ &  \nearrow \ \ \ \ \ \ \ \ \searrow & \ &  \nearrow \ \ \ \ \ \ \
\
\searrow
& \  & \nearrow \ \ \ \ \ \ \ \  \searrow & \ \cr
\ & \ & LP_k(\Psi)& \ & \Cal S_n(X, X\setminus Y) & \ & \ \cr
\ &  \searrow \ \ \ \ \ \ \ \ \nearrow & \ &  \searrow \ \ \ \ \ \ \
\
\nearrow
& \  & \searrow \ \ \ \ \ \ \ \  \nearrow & \ \cr
\rightarrow & LSP_{k} & \longrightarrow &
\Cal S_{l}(Y,Z,\eta) &
\longrightarrow & H_{l-1}(Y, \bold L_{\bullet}) & \rightarrow,
 \endmatrix
\tag 3.16
$$

$$
\matrix
\rightarrow & LS_{l}(F_Z) & \longrightarrow &
\Cal S_{n}(X,Y,Z) &
\rightarrow & \Cal S_{n}(X)& \rightarrow \cr
\ &  \nearrow \ \ \ \ \ \ \ \ \searrow & \ &  \nearrow \ \ \ \ \ \ \
\
\searrow
& \  & \nearrow \ \ \ \ \ \ \ \  \searrow & \ \cr
\ & \ & LSP_k& \ & \Cal S_n(X,Z,\nu) & \ & \ \cr
\ &  \searrow \ \ \ \ \ \ \ \ \nearrow & \ &  \searrow \ \ \ \ \ \ \
\
\nearrow
& \  & \searrow \ \ \ \ \ \ \ \  \nearrow & \ \cr
\rightarrow & \Cal S_{n+1}(X) & \longrightarrow &
LS_{k}(\Phi) &
\longrightarrow & LS_{l-1}(F_Z) & \rightarrow,
 \endmatrix
\tag 3.17
$$

and

$$
\matrix
\rightarrow & LS_{l+1}(F) & \longrightarrow &
LS_k(\Psi) &
\rightarrow & \Cal S_{n}(X,Y,Z)& \rightarrow \cr
\ &  \nearrow \ \ \ \ \ \ \ \ \searrow & \ &  \nearrow \ \ \ \ \ \ \
\
\searrow
& \  & \nearrow \ \ \ \ \ \ \ \  \searrow & \ \cr
\ & \ & \Cal S_{n+1}(X,Y, \xi)& \ & LSP_{k} & \ & \ \cr
\ &  \searrow \ \ \ \ \ \ \ \ \nearrow & \ &  \searrow \ \ \ \ \ \ \
\
\nearrow
& \  & \searrow \ \ \ \ \ \ \ \  \nearrow & \ \cr
\rightarrow & \Cal S_{n+1}(X,Y,Z) & \longrightarrow &
\Cal S_{n+1}(X) &
\longrightarrow & LS_{l}(F) & \rightarrow,
 \endmatrix
\tag 3.18
$$
where $l=n-q, k=n-q-q^{\prime}$.
The diagrams (3.15)--(3.18) are realized on the spectra level.
\endproclaim
\demo{Proof}
The transfer map  gives the commutative diagram (see \cite{18})
$$
\matrix
H_{n-q}(Y, \bold L_{\bullet})& \overset{\cong}\to{\to} &
H_{n}(X, X\setminus Y; \bold L_{\bullet})\\
                 &\searrow&\downarrow\\
& & H_{n-1}(X\setminus Y; \bold L_{\bullet}).\\
\endmatrix
\tag 3.19
$$
Consider the commutative triangle
$$
\matrix
LP_{n-q-q^{\prime}}(\Psi) &\to & L_n(\pi_1(X\setminus Y)\to \pi_1(X))\\
& \searrow  &\downarrow \\
 & & L_{n-1}(\pi_1(X\setminus Y)) \\
\endmatrix
\tag 3.20
$$
which lies in  the commutative diagram (3.4).

The results of \cite{18, Proposition 7.2.6} provide the maps
from the groups in diagram (3.19) to the corresponding groups
of diagram (3.20).
 On the spectra level
the cofibres of this map give a homotopy commutative
triangle  of spectra
$$
\matrix
\Bbb S(Y,Z, \eta) &\to & \Sigma^{-q}\Bbb S(X,X\setminus Y)\\
& \searrow  &\downarrow \\
 &  &\Sigma^{-q+1} \Bbb S(X\setminus Y). \\
\endmatrix
\tag 3.21
$$
  By \cite{20} the diagram (3.21) induces a homotopy commutative diagram
$$
\matrix
\Bbb S(Y,Z, \eta) &\to & \Sigma^{-q}\Bbb S(X,X\setminus Y)& \to &
\Sigma^{q^{\prime}+1} \Bbb LSP \\
\downarrow= & &\downarrow & &\downarrow \\
\Bbb S(Y,Z, \eta) &\to & \Sigma^{-q+1}\Bbb S(X\setminus Y)& \to &
\Sigma^{-q+1}\Bbb S(X,Y,Z)
\endmatrix
$$
in which the rows are cofibrations, and
 the right square is  a  pull-back.  The
homotopy long exact sequences of the maps of this square give
 the
braid (3.15). In a similar way,  the
maps from (3.19) to (3.20) provide the pull-back square
$$
\matrix
\Sigma^{q^{\prime}}\Bbb LP(\Psi)&\to &
\Sigma^{-q}\Bbb L(\pi_1(X\setminus Y)\to \pi_1(X))\\
\downarrow & &\downarrow \\
\Bbb S(Y,Z, \eta) &\to & \Sigma^{-q}\Bbb S(X,X\setminus Y)\\
\endmatrix
$$
in which the  cofibers of the vertical maps are homotopy equivalent to the
spectrum $Y_+\land \bold L_{\bullet}$. From this square we obtain
the braid of
exact sequences (3.16). The diagram (3.17) is obtained in
a similar way if we consider on the spectra level the homotopy
commutative triangle of the cofibers of the map from $H_n(X, \bold
L_{\bullet})$ to the triangle of natural forgetful maps
$$
\matrix
LT_{n-q-q^{\prime}} &\to & LP_{n-q-q^{\prime}}(\Phi)\\
& \searrow  &\downarrow \\
 & & L_{n}(\pi_1(X)) \\
\endmatrix
\tag 3.22
$$
which are obtained from square (3.10).
We obtain diagram (3.18) in a  way similar to  that of diagram (3.17).
 To do this we have to  consider the commutative triangle
$$
\matrix
LT_{n-q-q^{\prime}} &\to & LP_{n-q}(F)\\
& \searrow  &\downarrow \\
 & & L_{n}(\pi_1(X)) \\
\endmatrix
$$
instead of the triangle (3.22).
\qed
\enddemo
\smallskip

Let $Y^{n-q}\subset X^n$ be a manifold pair with $n-q\geq 5$ and $q\geq 3$.
Then by \cite{18} we have isomorphisms
$$
LS_n(F)\cong L_n(\pi_1(Y)), \ LP_n(F)\cong
L_{n+q}(\pi_1(X))\oplus L_n(\pi_1(Y)).
\tag 3.23
$$
Consider  the triple of manifolds (1.8) with the conditions
$$
n-q-q^{\prime}\geq 5, \ q\geq 3, \  q^{\prime}\geq 3.
\tag 3.24
$$
By \cite{14, Theorem 3} we have isomorphisms
$$
LT_{n-q-q^{\prime}}\cong L_{n}(\pi_1(X))\oplus L_{n-q}(\pi_1(Y))\oplus
L_{n-q-q^{\prime}}(\pi_1(Z)).
\tag 3.25
$$

Next
we obtain  similar results for the $LSP_*$@-groups.
\bigskip

\proclaim{Theorem 3.6} Suppose the  triple of manifolds (1.8) satisfy
the conditions (3.24). Then
$$
LSP_{n-q-q^{\prime}}(X,Y,Z)\cong  L_{n-q}(\pi_1(Y))\oplus
L_{n-q-q^{\prime}}(\pi_1(Z)).
$$
\endproclaim

\demo{Proof} The result follows by considering  the diagram
(3.9) and using the isomorphisms (3.25) and (3.23).
\qed
\enddemo
\bigskip

\proclaim{Theorem 3.7} Suppose the triple of manifolds (1.8) satisfy
the conditions $n-q-q^{\prime}\geq 5$ and  $q\geq 3$. Then
$$
LSP_{n-q-q^{\prime}}(X,Y,Z)\cong  LP_{n-q-q^{\prime}}(\Psi).
$$
\endproclaim

\demo{Proof} We have  isomorphisms
$$
LS_{n}(F_Z)\cong L_n(\pi_1(Y\setminus Z)), \
LS_{n}(F)\cong L_n(\pi_1(Y)), \
LS_{n}(\Phi)\cong L_n(\pi_1(Z)), \
$$
since $q\geq 3$.
The isomorphism
$$
LNS_n\cong L_{n+q^{\prime}}(\pi_1(Y\setminus Z)\to \pi_1(Y))
$$
follows from  diagram (2.18), since $q\geq 3$. The assertion of the
theorem follows now by chasing diagram
(3.8).
\qed
\enddemo

The Theorems 3.6 and 3.7 explain the geometrical
meaning of the obstruction groups $LSP_*$. These groups provide
 obstructions to  surgery on the submanifold pair
$(Y,Z)$ inside the ambient manifold $X$.

\subhead 4. Examples and applications
\endsubhead
\bigskip

A pair   of manifolds $Y\subset X$ is called a {\it Browder-Livesay pair}
if $Y$ is an one-sided submanifold of codimension 1 and the horizontal
maps in the square (1.3) are isomorphisms (see \cite{3}, \cite{5}, \cite{6},
\cite{11}, and \cite{12}). In this case the splitting obstruction groups
 are denoted by
$$
 LN_{n}(\pi_1(X\setminus
Y)\to \pi_1(X))=LS_{n}(F).
$$
Suppose the pairs of manifolds $(X,Y)$
and  $(Y,Z)$  in the triple (1.8)  are  Browder-Livesay
pairs. In this case $q=q^{\prime}=1$.
Denote by $r_p$ the map
$$
L_n(\pi_1(X))\to LSP_{n-3}(X,Y,Z)
$$
in the braid (3.4).
Let
$$
r:L_n(\pi_1(X))\to LS_{n-2}(F)=LN_{n-2}(\pi_1(X\setminus
Y)\to \pi_1(X))
$$
denote the map in the braid (1.5). The map $r$ gives the Browder-Livesay
invariant of an element $x\in L_n(\pi_1(X))$. If $r(x)\ne 0$ then the element
$x$ is not realized by a normal map of closed manifolds \cite{5}.

In the paper \cite{6} the invariants $A$ and $B$
 were defined. The invariant
$A$ coincides with $r$, and the invariant $B$ is defined on the kernel
of the invariant $A$. The invariant $B$ is called the
second Browder-Livesay invariant \cite{11}.
It is proved in \cite{6} that if
$B(x)\ne  0$ then the element $x$ is not realized by a normal map
of closed manifolds.
\bigskip

\proclaim{Proposition 4.1} Suppose the
 pairs of manifolds $(X,Y)$ and $(Y,Z)$ are  Browder-Livesay pairs. Then
$r_p(x)\neq 0$ if and only if $A(x)\ne 0$ or $B(x)\ne 0$.

\endproclaim
\demo{Proof} Consider the exact sequence fitting into the diagram (3.7)
$$
\cdots\to LT_{n-2}(X,Y,Z)\to L_n(\pi_1(X))\overset{r_p}\to{\to}
LSP_{n-3}(X,Y,Z)\to \cdots
$$
The proposition follows from this  exact sequence  and
 \cite{15, Theorem 3}.
\qed \enddemo
\smallskip

\proclaim{Corollary 4.2} If $r_p(x)\ne 0$
 then the element $x\in L_n(\pi_1(X))$ is not
realized by a normal map of closed manifolds.
\endproclaim
\smallskip

Next we compute some  $LSP$@-groups.
Consider the triple
$$
(Z\subset Y\subset X)=(\Bbb R\Bbb P^{n}\subset
\Bbb R\Bbb P^{n+1} \subset \Bbb R\Bbb P^{n+2})
\tag 4.1
$$
of real projective spaces
with $n\geq 5$.
The orientation homomorphism
$$
w:\pi_1(\Bbb R\Bbb P^{k})=\Bbb Z/2 \to
\{\pm 1\}
$$ is
trivial for
$k$ odd and nontrivial for $k$ even.
We have the following table for surgery obstruction groups
(see \cite{12} and \cite{21})
\bigskip
$$
\matrix
              & n=0  & n=1& n=2 & n=3 \cr
   L_n(1)     & {\Bbb Z}                   & 0 & {\Bbb Z}/2 & 0 \cr
L_n({\Bbb Z/2^+})&{\Bbb Z}\oplus {\Bbb Z}& 0 &{\Bbb Z}/2&{\Bbb Z}/2 \cr
L_n({\Bbb Z/2^-})&{\Bbb Z}/2& 0 &{\Bbb Z}/2& 0 \cr
\endmatrix
$$
The superscript "$+$" denotes the trivial orientation
of the corresponding group  and the superscript "$-$" denotes
the nontrivial orientation. For the Browder-Livesay pairs in (4.1) we
 have
the squares of fundamental groups
$$
F^{\pm}=
\left(\matrix
1 &\to & 1 \cr
\downarrow && \downarrow \cr
\Bbb Z/2^{\mp}&\to &\Bbb Z/2^{\pm}
\endmatrix\right).
$$
Furthermore we  have the isomorphisms (see \cite{12,  page 15} and \cite{21})
$$
LS_n(F^+)=LN_n(1\to \Bbb Z/2^+)=BL_{n+1}(+)=L_{n+2}(1)
$$
and
$$
LS_n(F^-)=LN_n(1\to \Bbb Z/2^-)=BL_{n+1}(-)=L_{n}(1).
$$
We recall intermediate computations of the
groups $LP_*(F^{\pm})$ and $LT^*(X,Y,Z)$
from \cite{14}.
The computation of $LP_*$@-groups for a pair $Y\subset X$ use
the  braid of exact sequences (1.5) (see, also \cite{19}).
The  natural map that forgets the
 manifold  $X$
$$
LS_{n}(F^{\pm})\to L_{n}(\Bbb Z/2^{\mp})
$$
coincides with the map
$$
l_n:BL_n(\pm)\to L_{n-1}(\Bbb Z/2^{\mp})
$$
in \cite{12, page 35}.
Using this result and
chasing the diagram (1.5)  we obtain the computations  (see also \cite{13})
$$
\matrix
LP_n(F^+)=&\Bbb Z/2,& \Bbb Z/2,& \Bbb Z/2, &\Bbb Z;\\
LP_{n}(F^-)=& \Bbb Z,& \Bbb Z/2, & \Bbb Z/2,&  \Bbb Z/2\\
\endmatrix
$$
for $n = 0, \ 1,\ 2,\ 3$ (mod $4$), respectively.

 Using connections  between  these groups and the
$LT_*(X,Y,Z)$@-groups,
the following result  was obtained in \cite{14}.

\proclaim{Proposition 4.3}
Let $M^{n-k}$ be a closed simply connected topological
manifold.  For the triple of manifolds
$$
(Z^n\subset Y^{n+1}\subset X^{n+2})=
(M^{n-k}\times\Bbb R\Bbb P^{k}\subset
M^{n-k}\times \Bbb R\Bbb P^{k+1} \subset
M^{n-k}\times \Bbb R\Bbb P^{k+2})
$$
with  $n\geq 5$, we have the following results.

For $k$ odd, the groups $LT_n$ are isomorphic to
$$
 \Bbb Z\oplus\Bbb Z/2, \Bbb Z/2, \Bbb Z\oplus\Bbb Z/2, \Bbb Z/2
$$
for $n = 0,1,2,3$ (mod $4$), respectively.

For $k$ even,
 $LT_0\cong \Bbb Z/2 \oplus \Bbb Z/2$
 and $LT_1\cong \Bbb Z/2$.
The groups $LT_3$ and $LT_2$
 fit into an exact sequence
  $$
  0\to LT_3 \to \Bbb Z \to \Bbb Z\to LT_2
\to \Bbb Z/2\to 0.
 $$
 \qed
\endproclaim
\smallskip

We apply these results to  compute
 $LSP_*$@-groups in the situation above.

\proclaim{Theorem 4.4} Under assumptions of the Proposition 4.3
 we have the following results.

For $k$ odd, the groups $LSP_n$ are isomorphic to
$$
 \Bbb Z, \Bbb Z, \Bbb Z/2, \Bbb Z/2
$$
for $n = 0,1,2,3$ (mod $4$), respectively.

For $k$ even, we have isomorphisms
 $LSP_0\cong LSP_1\cong \Bbb Z/2$.
The groups $LSP_3$ and $LSP_2$
 fit into an exact sequence
  $$
  0\to LSP_3 \to \Bbb Z \to \Bbb Z\to LSP_2 \to 0.
 $$
\endproclaim

\demo{Proof} Consider the case when $k$ is odd.
From \cite{14} we conclude that all the maps
$LT_n\to LP_{n+1}(F^+)$ are  epimorphisms. Now it is easy
 to describe the maps
$LP_n(F^+)\to L_{n+1}(\Bbb Z/2^+)$ in diagram (1.5). For
$n=1\bmod 4$ and $n=2\bmod 4$
these  maps are  isomorphisms $\Bbb Z/2\to \Bbb Z/2$.
For $n=0\bmod 4$ the map is trivial since the  group  $L_1(\Bbb Z/2^+)$
is trivial.
The map
$$
\Bbb Z=LP_3(F^+)\to L_0(\Bbb Z/2^+)=\Bbb Z\oplus\Bbb Z
$$
is an inclusion on a direct summand.  The image of this map
coincides with
the image of the map $L_0(1)\to L_0(\Bbb Z/2^+)$ that is induced by the
inclusion
$1\to \Bbb Z/2^+$. This follows
from the commutative triangle
$$
\matrix
& \Bbb Z &\\
& ||     & \\
& LP_3(F^+) & \\
 & \cong\nearrow \ \ \ \ \searrow  & \\
 L_0(1)& \overset{mono}\to{\longrightarrow} & L_0(\Bbb Z/2^+)\\
  ||    &            & ||   \\
  \Bbb Z & &       \Bbb Z\oplus \Bbb Z\\
\endmatrix
$$
in diagram (1.5).
From  diagram (3.4) we obtain an  exact  sequence
$$
\cdots\to LT_n\overset{\tau}\to{\to} L_{n+2}(\Bbb Z/2^+)\to LSP_{n-1}\to LT_{n-1}\to\cdots
$$
The map $\tau$ is the  composition
$$
LT_n\to LP_{n+1}(F^+)\to L_{n+2}(\Bbb Z/2^+)
$$
of maps that we already know.
Now we can compute the map $\tau$.
It  is trivial for $n=3$,
an isomorphism $\Bbb Z/2\to \Bbb Z/2$ for $n=1$, an epimorphism
$\Bbb Z\oplus \Bbb Z/2\to \Bbb Z/2$ with kernel $\Bbb Z$
for $n=0$, and  a homomorphism
 $\Bbb Z\oplus \Bbb Z/2\to \Bbb Z\oplus \Bbb Z$
with kernel $\Bbb Z/2$ and cokernel $\Bbb Z$ for $n=2$.
The result
for $k$ odd follows now from the exact sequence (4.24). The case of $k$ even is obtained
in a similar way.
\qed
\enddemo

\newpage

Information about Authors:

\noindent
Yuri V.  Muranov

\noindent
Department of General and Theoretical Physics,
Vitebsk State University,
Moskovskii pr. 33,
210026  Vitebsk,
Belarus

\noindent
e-mail:  ymuranov\@mail.ru

\bigskip

\noindent
 Rolando Jimenez

\noindent
Instituto de Matematicas, UNAM,
Avenida Universidad S/N, Col. Lomas de Chamilpa,
62210 Cuernavaca, Morelos,
Mexico

\noindent
e-mail: rolando\@aluxe.matcuer.unam.mx
\bigskip

\noindent
Du\v san Repov\v s:

\noindent
Institute for Mathematics, Physics and Mechanics, University of
Ljubljana, Jadranska 19, Ljubljana, Slovenia

\noindent
email: dusan.repovs\@uni-lj.si
\enddocument
\bye